\documentclass[12pt]{amsart}%
\usepackage{amsfonts}
\usepackage{graphicx}
\usepackage{amsmath}
\usepackage{amssymb}%
\input{amssym.def}
\newtheorem{theorem}{Theorem}[section]

\theoremstyle{remark}

\input epsf
\begin{document}
\title[Digons and angular derivatives]{Digons and angular derivatives of analytic self-maps of the unit disk}
\author[M.~D.~Contreras, S.~D{\'\i}az-Madrigal, A.Vasil'ev]{Manuel D.~Contreras, Santiago D{\'{\i}}az-Madrigal, \\ and Alexander Vasil'ev}

\address{M.~D.~Contreras, S.~D{\'\i}az-Madrigal: Camino de los Descubrimientos, s/n, Departamento de Matem\'atica Aplicada II, Escuela Superior de Ingenieros, Universidad de Sevilla, 41092, Sevilla, Spain}
\email{contreras@us.es, madrigal@us.es}

\address{A.Vasil'ev: Department of Mathematics, University of Bergen, Johannes Brunsgate 12, Bergen 5008, Norway}
\email{alexander.vasiliev@uib.no}

\thanks{Partially supported by Ministerio de Ciencia y Tecnolog{\'\i}a and the
European Union (FEDER) project BFM2003-07294-C02-02, by La
Consejer{\'\i}a de Educaci\'on y Ciencia de la Junta de
Andaluc{\'\i}a, and by the grants of the Norwegian Research
Council \#~177355, and of the University of Bergen}

\begin{abstract}
We present a geometric approach to a well-known sharp inequality, due to Cowen
and Pommerenke, about angular derivatives of general univalent self-maps of
the unit disk.

\end{abstract}
\maketitle

\section{Introduction}

Let $\mathbb{D}=\{z\in\mathbb{C}:\,\,|z|<1\}$ be the unit disk and
$\mathrm{Hol}(\mathbb{D},\mathbb{D})$ be the family of those analytic maps
defined on \thinspace$\mathbb{D}$ and with values also in $\mathbb{D}$.
Obviously, $\mathrm{Hol}(\mathbb{D},\mathbb{D})$ forms a semigroup with
respect to the functional composition with the identity map as the neutral
element. In this iteration context, the study of fixed points has always
played a prominent role. We recall that a point $\xi\in\overline{\mathbb{D}}$
is said to be a fixed point of a function $\varphi\in\mathrm{Hol}%
(\mathbb{D},\mathbb{D})$ whether $\lim_{r\rightarrow1}\varphi(r\xi)=\xi.$
Besides, when $\xi\in\partial\mathbb{D}$, this is indeed equivalent to assert
that the angular limit $\angle\lim\limits_{z\rightarrow\xi}\varphi(z)=\xi,$
that is to say, $\lim\limits_{z\rightarrow\xi,\,\,z\in\Delta_{\xi}}%
\varphi(z)=\xi$ for any Stolz angle $\Delta_{\xi}$ centred on $\xi$. Those
points $\xi\in\partial\mathbb{D}$ are usually called boundary fixed points of
$\varphi$. It is a remarkable fact that we can always talk about the
multiplier of this boundary fixed point $\xi$ in the sense that the following
angular limit exists
\[
\varphi^{\prime}(\xi):=\angle\lim\limits_{z\rightarrow\xi}\frac{\varphi
(z)-\xi}{z-\xi}.
\]
In fact, $\varphi^{\prime}(\xi)\in(0,+\infty)\cup\{\infty\}$ and when
$\varphi^{\prime}(\xi)\neq\infty\,,$ we say that $\xi$ is a regular (boundary)
fixed point. Moreover, these regular fixed points can be attractive if
$\varphi^{\prime}(\xi)\in(0,1),$ neutral if $\varphi^{\prime}(\xi)=1$ or
repulsive if $\varphi^{\prime}(\xi)\in(1,+\infty).$ Those fixed points $\xi
\in\mathbb{D}$ with $\,\,|\varphi^{\prime}(\xi)|<1$ are also named attractive.

A classical result by Denjoy and Wolff states that for a holomorphic self-map
$\varphi$ of the unit disk $\mathbb{D}$ other than a (hyperbolic) rotation,
there exists a unique fixed point $a\in\overline{\mathbb{D}}$ such that the
sequence of iterates $(\varphi_{n}(z))$ converges locally uniformly on
$\mathbb{D}$ to $a$ as $n\rightarrow\infty$. This point $a$ is called the
Denjoy-Wolff point of $\varphi$ and it is also characterized as the only fixed
point of $\varphi$ satisfying $\varphi^{\prime}(a)\in\mathbb{D}$. In other
words, $a$ is the only attractive fixed point of $\varphi$ in the above
multiplier sense. It follows from the hyperbolic metric principle that, if
$\varphi$ is not the identity, there can be no other fixed points in
$\mathbb{D}$ except the Denjoy-Wolff point but, nevertheless, $\varphi$ can
have many other repulsive or non-regular boundary fixed points.

One of the most studied properties concerning fixed points has been and is to
understand the relationship, if there is any, between the different fixed
points and their multipliers of a general function $\varphi\in\mathrm{Hol}%
(\mathbb{D},\mathbb{D})$. At this respect, Cowen and Pommerenke proved the
following nice result.

\begin{theorem}
\cite{CowenPommerenke} Let $\varphi\in\mathrm{Hol}(\mathbb{D},\mathbb{D})$ be
univalent with an attractive Denjoy-Wolff point $a\in\partial\mathbb{D}$, and let
$\xi_{1},\dots,\xi_{n}$ be $n\,\ $different repulsive boundary fixed points of
$\varphi.$ Then,%
\[
\sum\limits_{k=1}^{n}\frac{1}{\log\varphi^{\prime}(\xi_{k})}\leq-\frac{1}%
{\log\varphi^{\prime}(a)}.
\]
Moreover, this inequality is sharp.
\end{theorem}

As far as we know, there are two proofs of the above result in the literature.
One of them is the original one due to Cowen and Pommerenke and relies on the
so-called Grunsky-type inequalities. Later on, Li \cite{Li} obtained an
alternative proof based on the theory of contractively contained spaces,
sometimes also called de Branges-Rovnyak spaces. Undoubtedly, both approaches
have a definite heavy analytic flavour. In this paper, we provide a strictly
geometric proof of this inequality which, apart from its own interest, allows
us to improve it in a certain sense. Moreover, this point of view also
suggests to state it in a "weighted" way. Namely, we will prove the following.

\begin{theorem}
Let $\varphi\in\mathrm{Hol}(\mathbb{D},\mathbb{D})$ be univalent with
an attractive Denjoy-Wolff point $a\in\partial\mathbb{D}$, and let $\xi_{1}%
,\dots,\xi_{n}$ be $n\,\ $different repulsive boundary fixed points of
$\varphi.$ Likewise, consider $n$ non-negative numbers $\alpha_{1},\dots
,\alpha_{n}$, such that $\sum\limits_{k=1}^{n}\alpha_{k}=1.$ Then,%
\[
\prod_{k=1}^{n}\varphi^{\prime}(\xi_{k})^{\alpha_{k}^{2}}\geq|\varphi^{\prime
}(a)|^{-1}.
\]
Moreover, this inequality is sharp.
\end{theorem}

We want to mention that this theorem can also be considered as the `boundary'
version of a recent result due to Anderson and Vasil'ev \cite{Anderson}. At
the same time, a question naturally appears when comparing the above Theorems
1.1 and 1.2 and it is to explain the exact relationship between them. To this
aim, we need to use a simple fact from Calculus: consider the function
\[
f:D\subset\mathbb{R}^{n}\rightarrow\mathbb{R}\text{ \qquad}\alpha\mapsto
f(\alpha)=\sum_{k=1}^{n}\alpha_{k}^{2}c_{k},
\]
where $D:=\{\alpha=(\alpha_{1},\dots,\alpha_{n})\in\mathbb{R}^{n}:\alpha
_{1},...,\alpha_{n}\geq0$ and $\sum\limits_{k=1}^{n}\alpha_{k}=1\}$ and
$c_{1},...,c_{n}$ is an arbitrary family of $n$ non-negative numbers. Then,
applying Lagrange multiplier Theorem, we have that $f$ attains its minimum at the
point $\alpha=\frac{1}{\sum_{k=1}^{n}c_{k}}(c_{1},...,c_{n})\in D.$

Now, doing some computations and setting
\[
(\alpha_{1},\dots,\alpha_{n}):=\frac{1}{\sum_{k=1}^{n}\log\varphi^{\prime}%
(\xi_{k})}(\log\varphi^{\prime}(\xi_{1}),...,\log\varphi^{\prime}(\xi_{n}))
\]
in Theorem 1.2, we see that we can recover Theorem 1.1. At the
same time, this argument also shows that, putting together all the
weight inequalities, there is one which is the best and this one
is the original inequality of Cowen and Pommerenke. Anyhow, it is
worth recalling that all the weighted inequalities are sharp.
Indeed, in sections three and four we present two different
manners of obtaining those functions giving equality in such
weighted inequalities. These two ways correspond to techniques
from fractional iteration (section three) and from the theory of
quadratic differentials (section four).

\section{Digons, moduli, and proof}

\subsection{Digons and their reduced moduli}

A digon in the complex plane is a triple $(D,a,b),$ where $D$ is a
hyperbolically simply connected domain with piecewise smooth boundary and
$a,b$ are two finite points belonging to $\partial D$ called its vertices. As
usual, the vertices $a,b$ will be sometimes omitted when speaking about
digons. The main number associated with a digon is its reduced modulus
$m(D,a,b)$ which is defined by means of the concept of modulus of a simply
connected domain with respect to a family of paths lying in the domain (see
also \cite[Sections 2.1 and 2.2]{VasBOOK}). Not always this reduced modulus
exists but this will be the case in our context. For more general conditions
which guarantee the existence of this reduced modulus, we refer the reader to
\cite{Sol7} and \cite[section 2.2.3]{VasBOOK}.

Anyhow, more than working with single digons in the complex plane, we will be
really interested in dealing with certain families of non-overlapping digons
in the unit disk. Namely, consider an admissible system of curves on
$\mathbb{D},$ that is, a finite number $\gamma:=(\gamma_{1},\dots,\gamma_{n})$
of simple arcs not homotopically trivial on $\mathbb{D}$ that are not freely
homotopic pairwise on $\mathbb{D},$ which start and finish at fixed points of
$\partial\mathbb{D}$, are not homotopic to a point in $\mathbb{D}$ and do not
intersect. In this context, a system of non-overlapping digons $(D_{1}%
,\dots,D_{n})$ on $\mathbb{D}$ is said to be of homotopic type $(\gamma
_{1},\dots,\gamma_{n})$ if $(\gamma_{1},\dots,\gamma_{n})$ is an admissible
system of curves on $\mathbb{D}$ and for any $k\in\{1,\dots,n\},$ any arc on
$\mathbb{D}$ connecting the two vertices of $D_{k}$ is homotopic (not freely)
to $\gamma_{k}$.

Likewise, we fix a height vector $\alpha=(\alpha_{1},\dots,\alpha_{n})$ and we
will also require the digons to satisfy the so-called conditions of
compatibility of angles and heights or just associated conditions with
$\alpha$, i.e.,
\[
\delta_{k}(a_{k})=\frac{\pi\alpha_{k}}{\sum_{j\in I_{a_{k}}}\alpha_{j}}%
,\qquad\delta_{k}(b_{k})=\frac{\pi\alpha_{k}}{\sum_{j\in I_{b_{k}}}\alpha_{j}%
},\qquad\text{ }k=1,\dots,n,
\]
where $I_{c}$ is the set of indices which refers to the digons $D_{j}$ with a
vertex at $c$. At the same time, $\delta_{k}(a_{k}):=\sup\,\Delta_{a_{k}}$ and
$\delta_{k}(b_{k})=:\sup\Delta_{b_{k}}$ are the inner angles of the digon
$D_{k},$ that is, $\Delta_{a_{k}}$ (resp. $\Delta_{b_{k}})$ is the opening of
a Stolz angle inscribed in $D_{k}$ at $a_{k}$ (resp. $b_{k})$.

\subsection{Proof of Theorem 1.2}

Let $a$ be an attractive Denjoy-Wolff boundary point in $\partial \mathbb D$ as in Introduction.
Consider an admissible curve system $\gamma$ in $\mathbb{D}$ with
$\gamma:=(\gamma_{1},...,\gamma_{n}),$ and where each $\gamma_{k}$ is a simple
arc starting at $a$ and finishing at the repulsive fixed point $\xi_{k}.$
According to \cite{Sol7}, we can find an extremal system of non-overlaping
digons $(\widehat{D_{1}},...,\widehat{D_{n}})$ of homotopy type $\gamma$ and
satisfying the compatibility conditions associated with the vector
$\alpha:=(\alpha_{1},...,\alpha_{n})$
given in Theorem 1.2. Since $\varphi$ is univalent and $a,\xi_{1},...,\xi_{n}$
are regular boundary fixed points of $\varphi,$ we see that $(\varphi
(\widehat{D_{1}}),...,\varphi(\widehat{D_{n}}))$ is again a system of
non-overlaping digons of homotopy type $\gamma$ and satisfying the
compatibility conditions associated with $\alpha.$ In particular, the reduced
modulus is well-defined for each $\widehat{D_{k}}$ and each $\varphi
(\widehat{D_{k}})$. Hence, by the extremal character of the first system
\cite[Section 2.7.1]{VasBOOK}, we have that
\[
\sum\limits_{k=1}^{n}\alpha_{k}^{2}m(\varphi(\widehat{D_{k}}),a,\xi_{k}%
)\geq\sum\limits_{k=1}^{n}\alpha_{k}^{2}m(\widehat{D_{k}},a,\xi_{k}).
\]
Now, applying the change-variable formula \cite[Theorem 2.2.2]{VasBOOK}, we
deduce that
\begin{align*}
\sum\limits_{k=1}^{n}\alpha_{k}^{2}m(\widehat{D_{k}},a,\xi_{k})  &  \leq
\sum\limits_{k=1}^{n}\alpha_{k}^{2}m(\widehat{D_{k}},a,\xi_{k})\\
&  +\sum\limits_{k=1}^{n}\alpha_{k}^{2}\left(  \frac{1}{\delta_{k}(a)}%
\log\varphi^{\prime}(a)+\frac{1}{\delta_{k}(\xi_{k})}\log\varphi^{\prime}%
(\xi_{k})\right)  .
\end{align*}
Since%
\[
\delta_{k}(a)=\frac{\pi\alpha_{k}}{\sum_{j=1}^{n}\alpha_{j}}=\pi\alpha
_{k}\text{ and }\delta_{k}(\xi_{k})=\frac{\pi\alpha_{k}}{\alpha_{k}}%
=\pi,\text{ }%
\]
we conclude that
\[
0\leq\pi\sum\limits_{k=1}^{n}\alpha_{k}^{2}\left(  \frac{1}{\pi\alpha_{k}}%
\log\varphi^{\prime}(a)+\frac{1}{\pi}\log\varphi^{\prime}(\xi_{k})\right)
=\log\varphi^{\prime}(a)+\alpha_{k}^{2}\log\varphi^{\prime}(\xi_{k}).
\]
The proof ends by taking exponentials in the above inequality.

\section{Attaining the equality by means of fractional iteration}

Let us recall that a \textit{(one-parameter) semigroup of analytic functions}
is any continuous homomorphism $\Phi:t\mapsto\Phi(t)=\varphi_{t}$ from the
additive semigroup of non-negative real numbers into the composition semigroup
of all analytic functions which map $\mathbb{D}$ into $\mathbb{D}$. That is,
$\Phi$ satisfies the following three conditions:

\begin{enumerate}
\item[a)] $\varphi_{0}$ is the identity in $\mathbb{D},$

\item[b)] $\varphi_{t+s}=\varphi_{t}\circ\varphi_{s},$ for all $t,s\geq0,$

\item[c)] $\varphi_{t}(z)$ tends to $z$ locally uniformly in $\mathbb{D}$ as
$t\rightarrow0$.
\end{enumerate}

It is well-known that the functions $\varphi_{t}$ are always univalent. If $a$
is the Denjoy-Wolff point of one of the functions $\varphi_{t_{0}}$, for some
$t_{0}>0,$ then $a$ is the Denjoy-Wolff point of all the functions of the
semigroups, that is, all the functions of the semigroup share the Denjoy-Wolff
point. Moreover, if a point $\xi\in\partial\mathbb{D}$ is a boundary fixed
point of $\varphi_{t_{0}}$, for some $t_{0}>0,$ then it is a boundary fixed
point of all $\varphi_{t}$ \cite{ContrDiazPom}.

An important characteristic of a semigroup $\Phi$ is given by its Koenigs
functions. We will recall this notion just when the Denjoy-Wolff point $a$ is
in the boundary of the unit disk. In this case, there exists a unique
univalent function $h:\mathbb{D}\rightarrow\mathbb{C}$ such that $h(0)=0$, and
for every $t\geq0$,
\[
h\circ\varphi_{t}(z)=h(z)+t\qquad(z\in\mathbb{D})
\]
(see, for example, \cite{Shoikhet}). The function $h$ is called the Koenigs
functions of $\Phi$. Let us write $\Omega=h(\mathbb{D})$, and let $\nu(\Omega)$ be
the supremum of all positive numbers $y$ such that there is $c\in\Omega$ with
\[
\{c+it:\,\,-y/2<t<y/2\}\subset\Omega.
\]
The quantity $\nu(\Omega)$ is well defined and $\nu(\Omega)\in(0,\infty]$. We
have that $\nu(\Omega)<\infty$ if and only if $\varphi_{t}^{\prime}(a)<1$ for
some $t$ (if and only if $\varphi_{t}^{\prime}(a)<1 $ for all $t).$ In this
case, we have that
\[
\varphi_{t}^{\prime}(a)=\exp\left(  -\frac{\pi t}{\nu(\Omega)}\right)  ,
\]
for all $t\geq0$ \cite[Theorem 2.1]{ContrDiaz}. Let us denote $V(\Omega
)=\cap_{t\geq0}(\Omega+t)$. The set $V(\Omega)$ is called the invariant set
associated with $\Omega$. The invariant set can be empty. If it is not empty,
then it can be represented as at most a countable family of disjoint connected
components $V_{k}$. Each component $V_{k}$ is either an open strip or a
halfplane parallel to the real axis. We denote by $\beta_{k}(\Omega)$ the
width of the invariant strip $V_{k}$. By \cite[Theorems 2.5 and 2.6]%
{ContrDiaz}, there is a bijection from the set of all the invariant strips of
the invariant set and the set of repulsive regular boundary fixed points of
the functions of the semigroup. In fact, if $\xi_{k}$ is a repulsive regular
boundary fixed point of $\Phi$ with associated strip $V_{k}$, then
\[
\varphi_{t}^{\prime}(\xi_{k})=\exp\left(  \frac{\pi t}{\beta_{k}(\Omega
)}\right)  ,
\]
for all $t\geq0$.

Now, take positive numbers $\alpha_{1},\dots,\alpha_{n}$ with $\sum
\limits_{k=1}^{n}\alpha_{k}=1.$ Consider the domain
\[
\Omega =\{z\in \mathbb{C}:|\operatorname{Im}(z)|<\frac{1}{2}\}\setminus
\bigcup_{k=1}^{n-1}\{z\in \mathbb{C}:\operatorname{Im}(z)=-\frac{1}{2}%
+\sum_{j=1}^{k}\alpha _{j}\text{ and }\operatorname{Re}(z)<\gamma _{k}\},
\]
where $\gamma_{k}$ is a negative real number. Take $h$ a Riemann map of
$\Omega$ such that $h(0)=0$ and the semigroup given by $\varphi_{t}%
(z)=h^{-1}(h(z)+1)$ for all $t\geq0$ and $z\in\mathbb{D}.$ Then the
Denjoy-Wolff point of the functions of the semigroup is a point $a\in
\partial\mathbb{D}$, it has $n$ different repulsive regular boundary fixed
points $\xi_{1},\dots,\xi_{n},$ and
\[
\varphi_{t}^{\prime}(a)=\exp\left(  -\pi t\right)  \quad\text{and}\quad\text{
}\varphi_{t}^{\prime}(\xi_{k})=\exp\left(  \frac{\pi t}{\alpha_{k}}\right)
\quad\text{for all}\quad t>0.
\]
Then,
\[
\prod_{k=1}^{n}\varphi_t^{\prime}(\xi_{k})^{\alpha_{k}^{2}}=\exp\left(
{\textstyle\sum_{k=1}^{n}} \alpha_{k}^{2}\frac{\pi
t}{\alpha_{k}}\right)  =\exp\left(  \pi t\right)
=\frac{1}{\varphi_t^{\prime}(a)}.
\]
This proves that the inequality in Theorem 1.2 is sharp.

\section{Attaining the equality by means of quadratic differentials}

Coming back again to the proof of Theorem 1.2, we see that any $\varphi
\in\mathrm{Hol}(\mathbb{D},\mathbb{D})$ attaining the inequality given in the
theorem necessarily satisfies that
\[
\sum\limits_{k=1}^{n}\alpha_{k}^{2}m(\varphi(\widehat{D_{k}}),a,\xi_{k}%
)=\sum\limits_{k=1}^{n}\alpha_{k}^{2}m(\widehat{D_{k}},a,\xi_{k}),
\]
where $(\widehat{D_{1}},...,\widehat{D_{n}})$ is the extremal system of
non-overlapping digons associated with the curve system $(\gamma
_{1},...,\gamma_{n})$ and the vector $\alpha$. According to \cite[Section
2.7]{VasBOOK}, and in these circumstances, each $\widehat{D_{k}}$ is then a
strip domain in the trajectory structure of a unique quadratic differential
$Q(z)dz^{2}$ and there is a conformal map $g_{k}(z)$, $z\in\widehat{D_{k}}$
that satisfies the differential equation
\[
\alpha_{k}^{2}\left(  \frac{g_{k}^{\prime}(z)}{g_{k}(z)}\right)  ^{2}=4\pi
^{2}Q(z)
\]
and maps $\widehat{D_{k}}$ onto the strip $\mathbb{C}\setminus\lbrack
0,\infty)$. Moreover,
\[
Q(z)dz^{2}=A\frac{\prod_{j=1}^{n-1}(z-e^{i\beta_{k}})^{2}}{(z-a)^{2}%
\prod_{k=1}^{n}(z-\xi_{k})^{2}}dz^{2}.
\]
The complex constant $A$ is defined by the following conditions: the unit
circle is a trajectory of the differential, the integration of a branch of
$\sqrt{Q(z)}dz$ along the orthogonal trajectory about 0 gives 1, and the integration of a branch of $\sqrt{Q(z)}dz$ along the orthogonal
trajectory about 0 within each $\widehat{D_{k}}$ gives the height $\alpha_{k}$. Each point $e^{i\beta_{k}}$ lies in the arc of the unit circle joining the
points $\xi_{k}$ and $\xi_{k+1}$. Since this extremal configuration is unique,
so is the corresponding extremal function when fixing the points
$e^{i\beta_{k}}$. Indeed, this function $w=\widehat{\varphi}(z)$ satisfies the
complex differential equation
\[
\frac{dw}{dz}=\frac{(w-a)\prod_{j=1}^{n-1}(z-e^{i\beta_{k}})\prod_{k=1}%
^{n}(w-\xi_{k})}{(z-a)\prod_{j=1}^{n-1}%
(w-e^{i\beta_{k}})\prod_{k=1}^{n}(z-\xi_{k})},
\]
and the domain $\widehat{\varphi}(\mathbb{D})$ is the unit disk $\mathbb{D}$
minus at most $n-1$ analytic arcs starting at the points $e^{i\beta_{k}}$
along the trajectories of the quadratic differential $Q(w)dw^{2}$. Their
length depends on concrete values of the components of the vector $\alpha$ and
on the angular derivative $\varphi^{\prime}(a)$. In particular, all of the
extremal functions are of the same type as those described in section three.

The connection between the two approaches given in this section and section
three is clear. Namely, given the vector $\alpha$, the orbits of points under
the extremal dynamics considered in section 3 coincide with the trajectories
of the quadratic differentials in this section.


\begin{thebibliography}{9}                                                                                                %


\bibitem {Anderson}J. M. Anderson and A. Vasil'ev, \textit{Lower Schwarz-Pick
estimates and angular derivatives}, arXiv:math.CV/0608531 v2.

\bibitem {ContrDiaz}M.~D.~Contreras and S.~D{\'{\i}}az-Madrigal,
\textit{Analytic flows in the unit disk: angular derivatives and boundary
fixed points}, Pacific J. Math. \textbf{222} (2005), 253--286.

\bibitem {ContrDiazPom}M.~D.~Contreras, S.~D{\'{\i}}az-Madrigal, and
Ch.~Pommerenke, \textit{Fixed points and boundary behaviour of the Koenigs
function}, Ann. Acad. Sci. Fenn. \textbf{29} (2004), 471--488.

\bibitem {CowenPommerenke}C.~C.~Cowen and Ch.~Pommerenke, \textit{Inequalities
for the angular derivative of an anlytic function in the unit disk}, J.~London
Math. Soc. (2) \textbf{26} (1982), 271--289.

\bibitem {Li}K. Y. Li, \textit{Inequalities for fixed points of holomorphic
functions}, Bull. London Math. Soc. \textbf{22} (1990), 446-452.

\bibitem {Pomm75}Ch.~Pommerenke, \textit{Univalent functions, with a chapter
on quadratic differentials by G.~Jensen}, Vandenhoeck \& Ruprecht,
G\"ottingen, 1975.

\bibitem {Shoikhet}D.~Shoikhet, \textit{Semigroups in geometrical function
theory}, Kluwer Academic Publishers, Dordrecht, 2001.

\bibitem {Sol7}A.~Yu.~Solynin, \textit{Modules and extremal metric problems},
Algebra i Analiz, \textbf{11} (1999), no. 1, 1--86; English transl.,
St.-Petersburg Math. J., \textbf{11} (2000), no. 1, 1--70.

\bibitem {VasBOOK}A. Vasil'ev, \textit{Moduli of families of curves for
conformal and quasiconformal mappings.} Lecture Notes in Mathematics, vol.
1788, Springer-Verlag, Berlin--New York, 2002.
\end{thebibliography}
\end{document}